\documentclass[12pt]{article}
\usepackage{amsmath,amssymb,amsfonts,amsthm}

\author{Victor Turchin
\thanks{Partially supported by the grants NSH-1972.2003.01, MK-451.2003.01}}
\title{Calculus of the first non-trivial 1-cocycle of the space of long knots}
\date{}

\newtheorem{theorem}{Theorem}
\newtheorem*{lemma}{Lemma}
\newtheorem{remark}{Remark}[section]

\newtheorem*{conjecture}{Conjecture}
\def\mod{\mathop{\rm mod}\nolimits}
\newcommand\R{{\mathbb R}}
\newcommand\Z{{\mathbb Z}}
\newcommand\Q{{\mathbb Q}}
\newcommand\K{{\mathcal K}}

\newcommand\ie{{ i.e. }}
\newcommand\cf{{ cf. }}

\def\nbox{\quad$\Box$}
\def\nboxm{\quad\Box}
\def\B{{\textit B}}
\def\S{{\textit S}}

\def\DIAGRAMMAone{\begin{picture}(36,9)
                 \put(0,0){\line(1,0){36}}
                 \put(6,0){\line(1,1){9}}
                 \put(15,9){\vector(1,-1){9}}
                 \put(21,9){\vector(-1,-1){9}}
                 \put(21,9){\line(1,-1){9}}
                 \put(24,3){\vector(-2,-1){6}}
                 \put(24,3){\line(2,-1){6}}
                \end{picture}
               }

\def\DIAGRAMMAtwo{\begin{picture}(30,9)
                 \put(0,0){\line(1,0){30}}
                 \put(6,0){\line(1,1){6}}
                 \put(12,6){\vector(1,-1){6}}
                 \put(18,6){\vector(-1,-1){6}}
                 \put(18,6){\line(1,-1){6}}
                 \put(12,-6){$\stackrel{2}{\mapsto}$}
                \end{picture}
                }

\def\DIAGRAMMAthree{\begin{picture}(28,9)
                 \put(0,0){\line(1,0){28}}
                 \put(14,8){\vector(-1,-1){8}}
                 \put(14,8){\line(1,-1){8}}
                 \put(6,0){\line(4,3){4}}
                 \put(10,3){\vector(4,-3){4}}
                 \put(10.5,-8){\begin{picture}(12,6)(-1.5,0)
                               \thinlines
                               \put(2,3){\vector(-1,1){3}}
                               \put(2,3){\vector(-1,-1){3}}
                               \put(2.2,2){$\mapsto$}
                               \put(-1.5,1.3){\scriptsize 2}
                               \put(-1.5,3.7){\scriptsize 1}
                              \end{picture}
                             }
                \end{picture}
               }
\def\DiagrPt{\begin{picture}(10,3)
               \put(0,0){\line(1,0){10}}
               \put(5,0){\circle*{1}}
               \put(3.35,-5){$\stackrel{1}{\mapsto}$}
             \end{picture}
             }

\def\DiagrVtwo{\begin{picture}(30,6)
               \put(0,0){\line(1,0){30}}
               \put(6,0){\line(1,1){6}}
               \put(12,6){\vector(1,-1){6}}
               \put(18,6){\vector(-1,-1){6}}
               \put(18,6){\line(1,-1){6}}
             \end{picture}
            }

\def\DiagrSI{\begin{picture}(16,4)
              \put(0,0){\line(1,0){16}}
              \put(8,4){\vector(-1,-1){4}}
              \put(8,4){\vector(1,-1){4}}
             \end{picture}
            }

\def\DiagrUC{\begin{picture}(16,4)
              \put(0,0){\line(1,0){16}}
              \put(8,4){\vector(-1,-1){4}}
              \put(8,4){\line(1,-1){4}}
             \end{picture}
            }

\def\DiagrOCminus{\begin{picture}(16,8)
              \put(0,0){\line(1,0){16}}
              \put(8,4){\line(-1,-1){4}}
              \put(8,4){\vector(1,-1){4}}
              \put(4.3,4){$-$}
             \end{picture}
            }

\def\DiagrL{\begin{picture}(30,9)
             \put(0,0){\line(1,0){30}}
             \put(6,0){\line(1,1){9}}
             \put(15,9){\vector(1,-1){9}}
             \put(15,3){\vector(-1,-1){3}}
             \put(15,3){\vector(1,-1){3}}
            \end{picture}
            }
\def\DiagrLnew{\begin{picture}(24,6)
                \put(0,0){\line(1,0){24}}
                \put(6,0){\line(1,1){6}}
                \put(12,6){\vector(1,-1){6}}
                \put(12,0){\circle*{1}}
                \put(10.35,-5){$\stackrel{2}{\mapsto}$}
               \end{picture}
               }

\def\loopOverCross{\begin{picture}(27,10)
                    \qbezier(2,0)(5,0)(13.5,1)
                    \qbezier(13.5,1)(18.5,2)(18.5,5)
                    \qbezier(18.5,5)(18.5,9.5) (13.5,9.5)
                    \qbezier(13.5,9.5)(9.7,9.1)(9.4,5)
                    \qbezier(9.4,5)(9.6,2.6) (12,2)
                    \qbezier(16.2,0.8)(21,0)(25,0)
                   \end{picture}
                  }
\def\loopUnderCross{\begin{picture}(27,10)
                      \qbezier(25,0)(22,0)(13.5,1)
                      \qbezier(13.5,1)(8.5,2)(8.5,5)
                      \qbezier(8.5,5)(8.5,9.5) (13.5,9.5)
                      \qbezier(13.5,9.5)(17.3,9.1)(17.6,5)
                      \qbezier(17.6,5)(17.4,2.6) (15,2)
                      \qbezier(10.8,0.8)(6,0)(2,0)
                    \end{picture}
                   }

\def\casp{\begin{picture}(27,10)
              \qbezier(2,0)(3,0)(6,1)
              \qbezier(6,1)(13.5,4)(13.5,9.5)
              \qbezier(13.5,9.5)(13.5,4)(21,1)
              \qbezier(21,1)(24,0)(25,0)
          \end{picture}
          }

\def\trefoil{\begin{picture}(37,10)
              \qbezier(0,0)(7,0)(9,1.5)
              \qbezier(9,1.5)(10.6,2.8)(11.7,2.8)
              \qbezier(14.3,2.8)(17,2.6)(20,1.4)
              \qbezier(20,1.4)(23,0)(24,0)
              \qbezier(24,0)(25,0)(25.8,1)
              \qbezier(26.6,2.2)(28,10)(20,10)
              \qbezier(20,10)(13,10)(13,2.8)
              \qbezier(13,2.8)(13.5,0)(16,0)
              \qbezier(16,0)(16.5,0)(19,1)
              \qbezier(21,1.8)(24,3)(26.3,1.6)
              \qbezier(26.3,1.6)(27.3,1)(31,0.4)
              \qbezier(31,0.4)(33,0.1)(35,0)
              \qbezier(35,0)(36,0)(37,0)
             \end{picture}
            }

\begin{document}
\maketitle \sloppy

\begin{abstract}
For the space of long knots in $\R^3$, Vassiliev's theory defines the so called {\it finite order cocycles}.
Zero degree cocycles are finite type knot invariants. The first non-trivial cocycle of positive dimension in the space of long knots has 
dimension one and order three. We apply Vassiliev's combinatorial formula, given in~\cite{V4}, and find the value $\mod 2$ 
of this cocycle on the 1-cycles that are obtained by dragging knots one along the other or by rotating around a fixed line.
\end{abstract}

\setcounter{section}{-1}
\setcounter{theorem}{-1}

\section{Introduction}\label{s0}
\subsection{Finite type invariants and finite type cocycles}
Nowadays, the finite type knot invariants became widely used. (These invariants are also often called Vassiliev invariants.) 
They form a  filtered space. The invariants {\it of type $\leq k$} are those invariants that, being extended to singular knots (knots with a finite number 
of transversal self-intersections), vanish on all knots with $k+1$ self-intersections~\cite{V1,BN,ChDL,K}.

Historically, this simple and natural definition is based on a quite complicated geometrical construction of Vassiliev. 
The idea of this construction is to consider the space of all smooth maps $\R^1\rightarrow\R^3$ (not only embeddings)
with a fixed behaviour at infinity\footnote{Since the results of this paper concern only the space of long knots, \ie the space of smooth embeddings
$\R^1\hookrightarrow\R^3$ with a fixed behaviour at infinity, all the arguments are given for this situation only.} --- the maps are assumed to 
coincide with a fixed linear embedding outside some compact set of the line. This space is an affine space $R^\omega$ of some infinite dimension
$\omega$. The knot space $\K\subset\R^\omega$ is included in it; the complement  $\Sigma =\R^\omega\setminus\K$ is called {\it discriminant} and consists of the maps 
with self-intersections and singularities (points of  degeneration of the differential).

By the Alexander duality, one has the following isomorphism of the reduced (co)homology groups:
$$
\tilde H^*(\K)\simeq\tilde H_{\omega-*-1}(\bar\Sigma), \eqno(*)
$$
where $\bar\Sigma$ is the one-point compactification of $\Sigma$. Strictly speaking expression $(*)$ does not have a precise mathematical 
sense because $\omega$ is infinity. Thus we have to consider cycles in the discriminant of infinite dimension (but finite codimension). To
give a precise mathematical sense to the expression ($*$), one should consider finite-dimensional approximations of the space of long knots,
\ie finite dimensional subspaces $\R^N\subset\R^\omega$ that are in general position with the discriminant $\Sigma$, \cf~\cite{V1}.
For all $N$ the homology groups of the discriminants $\R^N\cap\Sigma$ are computed by means of spectral sequences. These spectral sequences (taken with a 
necessary shift of dimensions) stabilize and define the {\it stable spectral sequence}. By increasing the dimension $N$ we can exhaust all the 
finite-codimensional cycles of the discriminant. But, probably, not all of them would stabilize --- some of them could be lost by always hiding in the unstable
domain of the associated spectral sequence, \cf~\cite{V1}. The cocycles of the knot space dual to the stabilizing cycles of the discriminant are called
{\it cocycles of finite order} (or {\it finite type}). The zero dimension cocycles of this type are exactly the finite type invariants. The order of a cocycle equals 
the index $i$ of the filtration term $\sigma_{i}$ (of the resolved discriminant $\sigma$):   
$$
\sigma_0\subset\sigma_1\subset\sigma_2\subset\ldots \eqno(**)
$$
in which the dual cycle is contained. Roughly speaking, the order of a cocycle shows the complexity of the strata of $\Sigma$
that form the dual cycle.

This work is devoted to the first non-trivial finite order cocycle of positive dimension. We denote it by $v_3^1$. 
This cocycle is one-dimensional and has order three. V.~A.~Vassiliev called it the Teiblum-Turchin cocycle.

\subsection{Cocycle $v_3^1$. Vassiliev's results}
For the first time this cocycle was discovered by D.~Teiblum and the author in 1995 by computer calculations.
We found then only its \lq\lq main part\rq\rq 
\, lying in  $\sigma_3\setminus\sigma_2$. The dimensional restrictions imply  
that all higher stable differentials of the spectral sequence associated to ($**$) preserve this cocycle, which implies
its existence.   V.~A.~Vassiliev generalized the construction of this cocycle for spaces of long knots in the space $\R^n$ 
of arbitrary dimension $n\geq 3$, \cf~\cite{V2,V3}, 
and then found its explicit description $\mod 2$, \cf~\cite{V4}. To formulate the result of the paper~\cite{V4} (applied to the case $n=3$), 
we need some notation.

Suppose that the fixed linear embedding $\R^1\hookrightarrow\R^3$ (with which the knots coincide outside some compact subset of the line ---
this compact subset depends on a knot)  is the map $t\mapsto (t,t,0)$.
Suppose also that the direction  $\partial_z=(0,0,1)$ is the direction \lq\lq upward\rq\rq. If for a pair of points  $A$, $B$ the vector 
$\overrightarrow{AB}$ has the same direction as the vector $\partial_z$, we say that point $B$ is above point $A$, and point $A$ is below $B$. 
We consider the vertical projection of knots to the plane $(x,y)$. The direction $\partial_x=(1,0)$ in the plane $(x,y)$ 
will be called the direction \lq\lq to the right\rq\rq. 

\begin{theorem}[Vassiliev~\cite{V4}]\label{t0}
The $\mod 2$ value of $v_3^1$ on a generic $1$-cycle in the space of long knots is equal to the number of points of this $1$-cycle
(taken with appropriate multiplicities) corresponding to the knots
$f:\R^1\hookrightarrow \R^3$ such that one of three following conditions holds: 

(i) There are five points $a<b<c<d<e$ in $\R^1$ such that  $f(a)$ is above $f(d)$,   $f(e)$ is above
$f(b)$ and $f(c)$;

(ii) There are four points $a<b<c<d$ in $\R^1$ such that $f(a)$ is above  $f(c)$, $f(b)$ is below $f(d)$, 
and the projection of the derivative  $f'(b)$ to the plane $(x,y)$ is directed to the right;

(iii) There are three points $a<b<c$ in $\R^1$ such that $f(a)$ is above  $f(b)$ but below $f(c)$, and the angle formed by the projections of 
$f'(a)$ and $f'(b)$ to the plane $(x,y)$ does not contain the direction \lq\lq to the right\rq\rq.

These points  are taken with appropriate multiplicities. This multiplicity is equal to the number of configurations of points in the line
that satisfy one of the conditions (i), (ii), (iii): points $a$ and $d$ are not uniquely defined by (i); points $a$ and $c$ are not uniquely 
defined by (ii); if condition (iii) is satisfied its multiplicity is always one. Conditions (i) and (iii) can be
satisfied simultaneously; the multiplicities in this case are added. $\nboxm$
\end{theorem}

\subsection{Diagrammic description of the cocycle $v_3^1$}
The aim of the paper~\cite{V4}  is to give a general method for finding explicit combinatorial expressions for finite type cohomology classes.
Formulas that can be obtained by this method generalize Polyak-Viro-Goussarov formulas  (used for the finite type invariants, \cf~\cite{PV1,GPV}), 
and are also described by diagrams. The three diagrams on Figure~\ref{fig1} correspond, respectively, to conditions (i), (ii) and (iii) of Theorem~\ref{t0}.

\unitlength=1.2mm
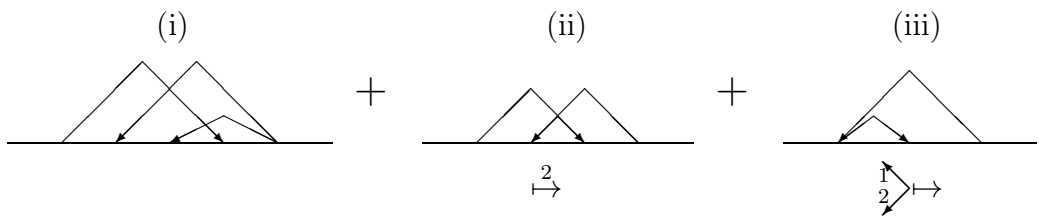
\begin{figure}[!h]
\begin{center}
 \begin{picture}(120,24)(0,-6)
   \put(16.5,12){(i)}\put(59.7,12){(ii)}\put(98,12){(iii)}
   \put(0,0){\begin{picture}(36,9)
              \put(0,0){\line(1,0){36}}
              \put(6,0){\line(1,1){9}}
              \put(15,9){\vector(1,-1){9}}
              \put(21,9){\vector(-1,-1){9}}
              \put(21,9){\line(1,-1){9}}
              \put(24,3){\vector(-2,-1){6}}
              \put(24,3){\line(2,-1){6}}
             \end{picture}
            }
   \put(38.5,4.5){\Large{+}}
   \put(46,0){\begin{picture}(30,9)
               \put(0,0){\line(1,0){30}}
               \put(6,0){\line(1,1){6}}
               \put(12,6){\vector(1,-1){6}}
               \put(18,6){\vector(-1,-1){6}}
               \put(18,6){\line(1,-1){6}}
               \put(12,-6){$\stackrel{2}{\mapsto}$}
              \end{picture}
             }
   \put(78.5,4.5){\Large{+}}
   \put(86,0){\begin{picture}(28,9)
               \put(0,0){\line(1,0){28}}
               \put(14,8){\vector(-1,-1){8}}
               \put(14,8){\line(1,-1){8}}
               \put(6,0){\line(4,3){4}}
               \put(10,3){\vector(4,-3){4}}
               \put(10.5,-8){\begin{picture}(12,6)(-1.5,0)
                           \thinlines
                           \put(2,3){\vector(-1,1){3}}
                           \put(2,3){\vector(-1,-1){3}}
                           \put(2.2,2){$\mapsto$}
                           \put(-1.5,1.3){\scriptsize 2}
                           \put(-1.5,3.7){\scriptsize 1}
                          \end{picture}
                         }
             \end{picture}
            }
 \end{picture}
\end{center}
\caption{cocycle $v_3^1$}\label{fig1}

\end{figure}

An arrow joining a couple of points means that these points have the same projection to the plane  $(x,y)$.
The direction of an arrow indicates which point is above and which is below: for example, in  diagram (i), the first point is above the fourth one,
and the fifth point is above the second and the third ones. The subscript $\stackrel{2}{\mapsto}$ under  diagram (ii) indicates that
the projection to the plane $(x,y)$ of the tangent vector at the second (from the left) point of the configuration is directed
\lq\lq to the right\rq\rq. The subscript
\begin{picture}(12,6)(-4,0)
                           \thinlines
                           \put(2,3){\vector(-1,1){3}}
                           \put(2,3){\vector(-1,-1){3}}
                           \put(2.2,2){$\mapsto$}
                           \put(-1.5,1.3){\scriptsize 2}
                           \put(-1.5,3.7){\scriptsize 1}
\end{picture}
means that the angle of the projections (to the plane $(x,y)$) of the derivatives in the first and second  points does not contain the direction
\lq\lq to the right\rq\rq. Further we will use several other specific notations (which will be explained later).

\section{Calculus of cocycle $v_3^1$}\label{s1}

\subsection{Primitivity}
Let $k$ be a long knot (which will be considered as a 0-cycle) and $c$ be a 1-cycle in the space of long knots. 
The space of long knots is an $H$-space: the product is defined as the composition of knots. Consider the products
 $k\ast c$ and $c\ast k$ of the above cycles.

\begin{remark}\label{r11}
{\rm Cycles $k\ast c$ and $c\ast k$ are homologous. Indeed, the cycle $c$ can be shrink to a very small size and dragged along
knot $k$.

A more general statement holds. In the paper~\cite{B}, R.~Budney has constructed an action of the operad of little discs on a space
weakly homotopy equivalent to the space of long knots. This action implies immediately a week homotopy commutativity of the 
multiplication on the space of long knots. \nbox }
\end{remark}

Note that none of the three diagrams from Figure~\ref{fig1} can be divided into two diagrams. Therefore the value of cocycle  
$v_3^1$ on cycles  $k\ast c$ and $c\ast k$ is equal to the value on cocycle $c$. In other words, the cocycle $v_{3}^1$ is primitive:
$$
\Delta v_3^1=v_3^1\otimes 1+1\otimes v_3^1,
$$
where $1$ denotes the unity of the cohomology algebra, \ie the invariant taking value $1$ on any knot.

\subsection{1-cycles}
In this section we describe some natural 1-cycles and evaluate  $v_3^1$ on these 1-cycles.

We consider cycles  of three types:

\medskip

1) Let $k$ be a long knot. Denote by $\hat k$ a 1-cycle that is obtained by rotation of $k$ around the axis $\{ (t,t,0)\,
|\, t\in\R\}$  (around the fixed straight line).

\bigskip

Recall that a {\it framing} of a long knot is a trivialization of its normal vector bundle that coincides 
with a fixed {\it constant trivialization} outside of some compact subset of the line (where the knot coincides already with the fixed linear
embedding). The normal vector bundle is oriented and twodimensional, hence the framing can be fixed by choosing only the first vector 
of the trivialization. The  {\it framing number} is defined as the linking number between the initial knot and the knot obtained by a slight displacement
of each point of the initial knot in the direction of the first trivialization vector. One says that a framing is  {\it trivial} if the framing number
equals zero.

\smallskip

2)  Let $k_1$ and $k_2$ be two knots. Let their composition $k_1\ast k_2$ be
an initial point of our 1-cycle. Fix a trivial framing on  $k_1$; then shrink knot $k_2$ and drag it to the left
along   $k_1$ as if knot $k_{2}$ were a small bead, and  turn this bead according to the framing of $k_{1}$. Now the 
knot $k_{2}$ gets stretched out to its initial size,  knots $k_{2}$ and $k_{1}$ get moved to the right ($k_{1}$ takes place of $k_{2}$, and 
$k_{2}$ takes place of $k_{1}$). Then we do the same procedure, but now we shrink the knot $k_{1}$ and drag it along the knot $k_{2}$ 
(according to a trivial framing of $k_{2}$). Finally, we come back to the initial point  $k_1\ast k_2$. The obtained 1-cycle will be 
denoted by $[k_1,k_2]$.

\begin{remark}\label{r12}
{\rm
As we already mentioned, the space $\K$ of long knots is endowed with an action of the operad of little squares. 
Therefore, the homology algebra $H_*(\K)$ is a Gerstenhaber algebra, \ie a graded commutative algebra with a Lie bracket $[\, .\, ,\, .\, ]$ 
that augments the degree by one and is compatible with the multiplication:
$$
[a,b\cdot c]=[a,b]\cdot c+ (-1)^{\tilde b\tilde c}[a,c]\cdot b.
$$
Our definition of the cycle $[k_1,k_2]$ is a particular case of this operation applied to the dimension zero homology.

The results of~\cite{B} imply that the Gerstenhaber algebra  $H_*(\K,\Q)$ is a free (Gerstenhaber) algebra, the space of its generators is  
the direct sum of the homology groups  $H_*(\K_p,\Q)$ of the space  $\K_p\subset\K$ of prime knots 
($K_{p}$ is the disjoint union of the connected components of $\K$ that correspond to prime knots), see also~\cite{BC}.
\nbox }
\end{remark}

\bigskip

3) Let us apply the previous construction to two identical knots $k$ and $k$.  Note that even after the first dragging of one 
knot along the other we come back to the initial point $k\ast k$. The obtained cycle will be denoted by  $k\circ k$.
Obviously, $[k,k]=2k\circ k$.

\medskip

It is easy to see that the definition of  cycles  $\hat k$, $[k_1,k_2]$, and $k\circ k$ is functorial. Therefore, the value of any 1-cocycle on these cycles 
must be an invariant of a knot (in the first and the second situations) or a biinvariant,\ie invariant of a pair of knots (in the second situation).

We will prove the following theorem:

\begin{theorem}\label{t1}
The value of the cocycle  $v_3^1$ (reduced $\mod 2$) on the cycles $[k_1,k_2]$ and $k\circ k$ is equal to zero, 
and on the cycle $\hat k$ is equal $\mod 2$ to the value of the first non-trivial Vassiliev invariant $v_2$ (of order $2$) on the knot
$k$. \nbox
\end{theorem}

Invariant $v_2$ is called the  {\it Casson invariant}, and its $\mod 2$ reduced version is called the {\it Arf invariant}, \cf~\cite{PV1,G,Kauf,L1,Ng,PV2}. 
The value $v_2(k)$ of this invariant on a knot can be expressed in terms of the plane projection of the knot according to the diagram:

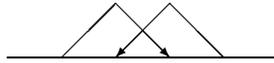
\begin{center}
\begin{picture}(30,6)
               \put(0,0){\line(1,0){30}}
               \put(6,0){\line(1,1){6}}
               \put(12,6){\vector(1,-1){6}}
               \put(18,6){\vector(-1,-1){6}}
               \put(18,6){\line(1,-1){6}}
\end{picture}
\end{center}

In other words,  $v_2(k)$ is equal to the number of configurations (counted with appropriate signs) of four points 
$a<b<c<d$ in $\R^1$ such that $f(a)$ is above $f(c)$, and $f(d)$ is above $f(b)$, where
$f:\R^1\hookrightarrow\R^3$ is the map that defines knot $k$.

\bigskip

It is natural to pose the following question: What would happen if, in the second and in the third constructions, one considers 
draggings along non-trivial framings\footnote{For example, V.~A.~Vassiliev  calculates in the paper~\cite{V4} the value of cocycle 
$v_3^1$ on a cycle obtained by dragging a trefoil knot along itself, but the framing number in the considered situation is not zero.}?
It is easy to see that in the first case one gets a cycle homologous to the sum:
$$
[k_1,k_2]+ N_1\cdot k_1\ast \hat k_2 +N_2\cdot k_2\ast\hat k_1,
$$
where $N_1$ and $N_2$ are, respectively, the framing numbers of the first and of the second knots. In the second case, the cycle is homologous 
to the sum: 
$$
k\circ k+N\cdot k\ast \hat k,
$$
where $N$ is the framing of $k$. By Theorem~\ref{t1} and the primitivity propriety of $v_{3}^1$, the value of $v_{3}^1$ on these 
cycles is equal  $\mod 2$ to  $N_1\cdot v_2(k_2)+N_2\cdot v_2(k_1)$ and $N\cdot
v_2(k)$, respectively.

\subsection{Flat representation of 1-cycles}
Let $k$ be a knot with a generic projection to the plane $(x,y)$. One can assign a {\it flat framing} 
of the knot  $k$ to this projection. The first trivialization vector of this framing is defined to be always parallel to the plane 
$(x,y)$ (in all the points of  knot $k$). Flat framings are more convenient for description of cycles  $[k_1,k_2]$ and
$k\circ k$. But usually these framings are not trivial. To make a flat framing trivial, we add a necessary number 
of loops  \unitlength=0.17em \loopOverCross
or \loopUnderCross to the plane projection of a knot. The first loop augments the flat framing number by one, and the second one diminishes this number by one. 
Representations of the cycles  $[k_1,k_2]$ and $k\circ k$ that use flat framings
will be called {\it flat representations}.

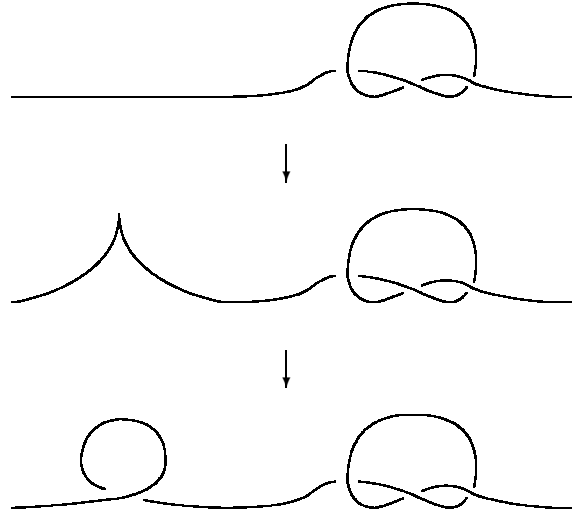
\begin{figure}[!ht]
\unitlength=0.3em
\begin{center}
\begin{picture}(60,10)
 \put(0,0){\line(1,0){23}}
 \put(23,0){\trefoil}
\end{picture}
\\
\vspace{5mm}
\begin{picture}(1,4)
 \put(0,3){\vector(0,-1){4}}
\end{picture}
\\
\vspace{4.5mm}
\begin{picture}(60,10)
 \put(-2,0){\casp}
 \put(23,0){\trefoil}
\end{picture}
\\
\vspace{5mm}
\begin{picture}(1,4)
 \put(0,3){\vector(0,-1){4}}
\end{picture}
\\
\vspace{4.5mm}
\begin{picture}(60,10)
 \put(-2,0){\loopOverCross}
 \put(23,0){\trefoil}
\end{picture}
\end{center}
\caption{twisting a loop}\label{fig2}
\end{figure}

Let us  also define a  {\it flat representation} of the cycle $\hat k$. Consider a knot $k$. We twist a loop on the left of $k$, see Figure~\ref{fig2}.

Then we shrink the knot and drag it to the left along the loop according to the loop's flat framing. After the knot  passes 
on the left, we stretch it back, untwist the loop and put the knot to its initial position. It is easy to see that the obtained cycle is homologous to the cycle
$\hat k$. It is also obvious that we do not encounter the configurations of Theorem~\ref{t0} while twisting or untwisting the loop. Therefore we need 
to compute the number of such configurations only while (flat) dragging the knot along the loop.

\subsection{Proof of Theorem~\ref{t1}}
Consider two knots $\B$ (big) and $\S$ (small) with generic projections to the plane 
$(x,y)$, that are consecutively tied on the same line. Let us compute the number of configurations mentioned in 
Theorem~\ref{t0} that we meet while dragging  knot 
 $\S$ along knot  $\B$ according the flat framing. 

Without loss of generality we may 
suppose that there are two possible  {\it standard} types of crossing points of knots' projections:

\unitlength=0.1em
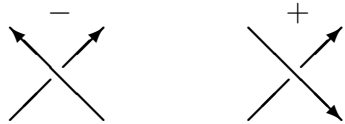
\begin{figure}[!h]
\thicklines
\begin{center}
\begin{picture}(30,35)
  \put(0,0){\line(1,1){13}}
  \put(17,17){\vector(1,1){13}}
  \put(30,0){\vector(-1,1){30}}
  \put(12,32){$-$}
\end{picture}
\begin{picture}(40,35)
\end{picture}
\begin{picture}(30,35)
 \put(0,0){\line(1,1){13}}
 \put(17,17){\vector(1,1){13}}
 \put(0,30){\vector(1,-1){30}}
 \put(12,32){$+$}
\end{picture}
\end{center}
\caption{standard crossings}\label{fig3}
\end{figure}

The direction of the derivative in  underpass points is  $(1,1)$ in the both negative  
$-$ and positive $+$ cases. The direction of the derivative in overpass points must be 
$(-1,1)$ in the negative $-$ case, and must be $(1,-1)$ in the positive $+$ case.

The 1-chain obtained by this dragging will be denoted  $\B\circ \S$. (If the flat framings of knots $\B$ and
$\S$ are trivial, then  $[\B,\S]=\B\circ\S+\S\circ\B$.)

The left column of Table~1 describes possible situations of appearance of the configurations described in Theorem~\ref{t0}.
A subscript $\S$ or $\B$ under a point of a diagram means to which knot (to small knot $\S$ or to big knot $\B$) this point belongs. 

\unitlength=1.18mm
\begin{table}[!h]
\caption{multiplication rule}\label{tt1} \vspace{3mm}
\begin{center}
\begin{tabular}{|cc|c|c|}
\hline
&$\B\circ\S$&$\B$&$\S$\\
\hline \raisebox{30pt}{\Large I} &
\begin{picture}(36,18)
 \put(0,5){\DIAGRAMMAone}
 \put(4.6,3.2){{\tiny $\B$}}
 \put(10.6,3.2){{\tiny $\S$}}
 \put(16.6,3.2){{\tiny $\S$}}
 \put(22.6,3.2){{\tiny $\B$}}
 \put(28.6,3.2){{\tiny $\B$}}
\end{picture}
&
\begin{picture}(30,11)
\put(0,5){\DiagrVtwo}
\end{picture}
&
\begin{picture}(16,11)
\put(0,5){\DiagrSI}
\end{picture}
\\
\hline
\raisebox{30pt}{\Large II} &
\begin{picture}(36,18)
 \put(0,5){\DIAGRAMMAone}
 \put(4.6,3.2){{\tiny $\S$}}
 \put(10.6,3.2){{\tiny $\S$}}
 \put(16.6,3.2){{\tiny $\S$}}
 \put(22.6,3.2){{\tiny $\S$}}
 \put(28.6,3.2){{\tiny $\B$}}
\end{picture}
&
\begin{picture}(16,11)
\put(0,5){\DiagrUC}
\end{picture}
&
\begin{picture}(30,16)
\put(0,5){\DiagrL}
\end{picture}
\\
\hline \raisebox{31pt}{\Large III} &
\begin{picture}(30,16)
 \put(0,6){\DIAGRAMMAtwo}
 \put(4.6,4.2){{\tiny $\S$}}
 \put(10.6,4.2){{\tiny $\S$}}
 \put(16.6,4.2){{\tiny $\S$}}
 \put(22.6,4.2){{\tiny $\S$}}
\end{picture}
&
\begin{picture}(10,12)
\put(0,6){\DiagrPt}
\end{picture}
&
\begin{picture}(30,12)
\put(0,6){\DiagrVtwo}
\end{picture}
\\
\hline \raisebox{31pt}{\Large IV} &
\begin{picture}(30,16)
 \put(0,6){\DIAGRAMMAtwo}
 \put(4.6,4.2){{\tiny $\S$}}
 \put(10.6,4.2){{\tiny $\S$}}
 \put(16.6,4.2){{\tiny $\S$}}
 \put(22.6,4.2){{\tiny $\B$}}
\end{picture}
&
\begin{picture}(16,12)
 \put(0,6){\DiagrUC}
\end{picture}
&
\begin{picture}(24,13)
 \put(0,6){\DiagrLnew}
\end{picture}
\\
\hline \raisebox{31pt}{\Large V} &
\begin{picture}(30,16)
 \put(0,6){\DIAGRAMMAtwo}
 \put(4.6,4.2){{\tiny $\B$}}
 \put(10.6,4.2){{\tiny $\S$}}
 \put(16.6,4.2){{\tiny $\B$}}
 \put(22.6,4.2){{\tiny $\B$}}
\end{picture}
&
\begin{picture}(30,12)
 \put(0,6){\DiagrVtwo}
\end{picture}
&
\begin{picture}(10,12)
\put(0,6){\DiagrPt}
\end{picture}
\\
\hline \raisebox{37pt}{\Large VI} &
\begin{picture}(28,19)
 \put(0,8){\DIAGRAMMAthree}
 \put(4.6,6.2){{\tiny $\S$}}
 \put(12.6,6.2){{\tiny $\S$}}
 \put(20.6,6.2){{\tiny $\B$}}
\end{picture}
&
\begin{picture}(16,15)
\put(0,8){\DiagrUC}
\end{picture}
&
\begin{picture}(16,15)
\put(0,8){\DiagrOCminus}
\end{picture}
\\
\hline
\end{tabular}
\end{center}
\end{table}

It is easy to write out the rules that exclude all the other situations:

1) The points that belong to knot $\S$ must be consecutive;

2) If one has an arrow joining a point from $\S$ with a point from $\B$ then all such arrows must have the same orientations: 
whether from  $\S$ to $\B$, whether from  $\B$ to $\S$ (This situation corresponds to moments when knot $\S$ passes through a crossing point
of $\B$. In the first case $\S$ passes through an overpass point, in the second case --- through an underpass point.);

3) In the case of a triple self-intersection precisely one point must belong to $\B$, and precisely two points must belong to $\S$ (diagrams (i) and (iii));

4) All the points of knot $\B$ that are joint to those of $\S$ must be consecutive;

5) If in some point $X$ of a diagram one has a restriction on tangent vector (projection directed to the right), then
this point must belong to  small knot  $\S$. If in this situation there are arrows joining $\S$ with $\B$, then the second point (that is joint with $X$)
must belong to  knot $\B$ (diagram (ii)).

\medskip

In Table~1 one can find two new notations.

In the first line \DiagrSI --- a two-sided arrow means that the points of a knot are in the same vertical line: the first point might be either 
over or below the second one.

In the line VI  --- a minus over a diagram means that one has a negative crossing point, see Figure~\ref{fig3}.

The assumption that all the crossing points are standard, see Figure~\ref{fig3}, implies the following lemma:

\begin{lemma} {\rm (multiplication rule)}
In all the six situations I-VI the number of appearances of configurations  from the left column in 1-chain $\B\circ \S$ is equal to the product of two numbers:
the first one being the number
of configurations  that correspond to the diagram from the same line, and the second column computed for knot  $\B$;  and the second one being the number 
of configurations that correspond to the diagram from the same line, but from the third column computed for knot  $\S$. \nbox
\end{lemma}

Consider, for example, the first line. It is easy to see that configuration  
 \unitlength=0,22em
 \begin{picture}(36,14)
 \put(0,1){\DIAGRAMMAone}
 \put(4.5,-1.6){{\tiny $\B$}}
 \put(10.5,-1.6){{\tiny $\S$}}
 \put(16.5,-1.6){{\tiny $\S$}}
 \put(22.5,-1.6){{\tiny $\B$}}
 \put(28.5,-1.6){{\tiny $\B$}}
\end{picture} in our chain $\B\circ\S$ appears only while  small knot $\S$ passes through the second point of configuration 
 \DiagrVtwo of  big knot $\B$ and exactly at the moment of its self-intersection \DiagrSI. So, the multiplication rule follows. The other cases
 are examined analogously.
 
 \medskip

To complete the proof of the theorem let us note that the following equalities hold:

\begin{gather}
II+IV+VI=0\mod 2.\notag\\
I+V=0\mod 2.\notag
\end{gather}

In both cases the common factor that corresponds to knot  $\B$ is taken out the parenthesis. And what is inside the parentheses
corresponds to the knots  $\S$, and is equal $\mod 2$ to zero. Only the summand  III remains. But the first factor in III is equal $\mod 2$ 
to the flat framing number of  $\B$, and the second factor --- to the value of the Arf invariant (simplest Vassiliev $\Z_{2}$-invariant of order 2)
on knot  $\S$.

\subsection{Open questions}
The first natural question is to find a similar explicit description of cocycle $v_3^1$ over $\Z$ ant to evaluate this cocycle on cycles 
 $[k_1,k_2]$, $k\circ k$ and $\hat k$.

\begin{conjecture}
The value of the cocycle $v_3^1$ (over $\Z$) is equal to zero on cycles $[k_1,k_2]$, $k\circ k$, and is equal to $v_{2}(k)$ (the value of the Casson invariant 
on knot $k$) on cycle $\hat k$. \nbox
\end{conjecture}

\medskip

The main propriety of cycles  $[k_1,k_2]$, $k\circ k$, $\hat k$ is the functoriality. In other words, the cycles do not depend
on the representatives of the knot isotopy classes. There exists a one more type of 1-cycles that has the same propriety --- 
we will denote it by  $\tilde k$. This 1-cycle is obtained from a knot  $k$ by the following way: Let us fix a trivial framing on knot $k$.
Consider a one-point compactification  $S^3$ of the space  $\R^3$. The space of long knots can be regarded as a space of knots in
$S^3$ (embeddings $S^1\hookrightarrow S^3$) with  fixed initial point and fixed derivative in the initial point. In the space of such knots one can 
consider a 1-cycle that is obtained by dragging knot $k$ through the initial point as if it was of wire that can not be deformed, \cf~\cite{H1}. 
The rotation angle at each moment is determined by the framing. On the plane projection it looks like the knot is turned inside out, \cf~\cite{H1}. 
The natural question is whether the value of the cocycle  $v_3^1$ on these 1-cycles is a finite type knot invariant.

\section*{Acknowledgement}
I would like to thank V.~A.~Vassiliev and S.~S.~Anisov for interesting discussions and for correcting the earlier versions of this paper.

\pagebreak

\vspace{8mm}

\noindent Independent University of Moscow\\
Moscow 121002, B. Vlassievskij per. 11\\
\tt  turchin@mccme.ru

\end{document}